\documentstyle[12pt]{article}
\def\vs{\vspace*}
\def\qed{\hfill$\Box$}

    \def\dim{{\rm dim}} 
\def\F{\hbox{$I\hskip -4pt F$}}
\def\dis{\displaystyle}
\def \Z{\hbox{$Z\hskip -5.2pt Z$}}
\def\sZ{\hbox{$\sc Z\hskip -4.2pt Z$}}
\def\ZZ{{\Z}\!\times\!{\Z}} 
\def\vsp{\vspace*{4pt}}
\def\Vir{\mbox{Vir}}

  \def\0{{\{0\}}}
\def\oz{\oplus_{i\in \sZ}}  \def\a{\alpha}

\def\AA{{\cal A}}
\def\BB{{\cal B}}
\def\CC{{\cal C}}
\def\EE{{\cal D}}
\def\nl{\newline}
\def\b{\beta}   
 
\def\sc{\scriptstyle}
\def\ssc{\scriptscriptstyle}
\def\bs{\backslash}
\def\cl{\centerline}
\begin{document}
%
%
\cl{{\bf $\ZZ$-graded Lie algebras generated by the Virasoro
 algebra and $sl_2\,$}\footnote{AMS Subject Classification:
 17B36, 17B65, 17B70\vs{4pt}\nl\hspace*{4.5ex}Supported by
a NSF grant 10171064 of China and a Fund from National Education
Ministry of China.}} \vs{4pt}\cl{(Appeared in {\it Math.~Nach.} 
{\bf 246/247} (2002), 188-201.)}
 \vskip 9pt
\par
\cl{
  Yucai Su$^{\dag}$,
  Jiangfeng Zhang$^{\dag}$
and
 Kaiming Zhao$^{\ddag}$
}\par
{\small
 {\it
  $^{\dag}$Department of Mathematics, Shanghai
  Jiaotong University, Shanghai 200030, P.~R. China
  \vskip -3pt\par
  $^{\ddag}$Institute of Mathematics,  Academy of Mathematics
  and System Sciences, Chinese Academy of Sciences,  Beijing 100080,
  P.~R.~China
}
}
\vskip 9pt
{\small
{\bf ABSTRACT} The present paper is another step toward the classification
of $\ZZ$-graded Lie algebras:
we classify  $\ZZ$-graded Lie algebras
$\AA=\oplus_{i,j\in\sZ}\AA_{i,j}$ over a field $\F$ of characteristic 0 with
$\dim\AA_{i,j}\le1$ for $i,j\in\Z$, satisfying:
\nl\hspace*{3ex}
(1) $\AA_0=\oplus_{i\in \sZ}\AA_{i,0}\simeq\Vir$, the centerless Virasoro
 algebra,
\nl\hspace*{3ex}
(2) $\AA_{0,-1},\AA_{0,0},\AA_{0,1}$ span the 3-dimensional simple Lie
 algebra $sl_2$,
\nl\hspace*{3ex}
(3) $\AA$ is generated by the centerless Virasoro algebra and $sl_2$.
\vspace*{-2pt}\par
These algebras include some rank 2 centerless Virasoro algebras,
and some rank 2 Block algebras and their extensions.
}
\par\vskip 10pt
\cl{{\bf\S1.  Introduction}}
\par
The classification of simple $\Z$-graded Lie algebras of
finite growth was given by O. Mathieu [M] in 1992.
It is natural to ask what would happen for $\ZZ$-graded Lie algebras.
This problem was initiated and studied by Osborn and Zhao. They
investigated several cases in [OZ1-OZ5] under some mild restrictions.
There are still some cases left untouched.
The present paper is another step by the authors toward
the classification of $\ZZ$-graded Lie algebras.
\par
Throughout this paper we assume that the $\ZZ$-graded Lie algebra
$\AA=\oplus_{i,j\in\sZ}\AA_{i,j}$ over a field $\F$ of characteristic 0 with
$\dim\AA_{i,j}\le1$ for $i,j\in\Z$, satisfies the following three conditions:
\nl\hspace*{3ex}
(1) $\AA_0=\oplus_{i\in\sZ}\AA_{i,0}\simeq\Vir$, the centerless Virasoro
 algebra;
\nl\hspace*{3ex}
(2) $\AA_{0,-1},\AA_{0,0},\AA_{0,1}$ span the 3-dimensional simple Lie
 algebra $sl_2$;
\nl\hspace*{3ex}
(3) $\AA$ is generated by $\AA_{0,\pm1}$ and $\AA_{0}$.
\par
To state our main theorem, we recall some known Lie algebras and construct
some other Lie algebras from them. For any $\a\in\F$ with $\a\ne0$,
we have the generalized centerless Virasoro algebra
$\Vir(\a)=\oplus_{i,j\in\sZ}L_{i,j}$ with bracket
$$[L_{i,j},L_{k,\ell}]=(k-i+(\ell-j)\a)L_{i+k,j+\ell}\mbox{ for }
(i,j),(k,\ell)\in \ZZ.
\eqno(1.1)$$
\par
For any $\a,\b\in\F$ with $\a\b\ne0$, denote
$Z=\ZZ\bs\{(-\a,\b),(-2\a,2\b)\}.$
We have (see [B]) the Block algebra $\BB(\a,\b)=\oplus_{(i,j)\in Z}\F L_{i,j}$ with the
bracket
$$
[L_{i,j},L_{k,\ell}]=
\left|\matrix{i+\a&j-\b\vspace*{2pt}\cr k+\a&\ell-\b\cr}\right|
L_{i+k,j+\ell}\mbox{ for }(i,j),(k,\ell)\in Z.
\eqno(1.2)$$
{}From [DZ], we see that
the Lie algebra $\BB(\a,\b)$ has the following possible central
extensions $\BB(\a,\b;a_1,a_2,a'_2)=\BB(\a,\b)\oplus \F{\bf c}_1
\oplus \F{\bf c}_2$
for any $a_1,a_2,a'_2\in \F$:
$$
[L_{i,j},L_{k,\ell}]=
\left|\matrix{i+\a&j-\b\vspace*{2pt}\cr k+\a&\ell-\b\cr}\right|
L_{i+k,j+\ell}
+\delta_{i+k,-\a}\delta_{j+\ell,\b}(\a j+\b i)a_1{\bf c}_1\hskip 3cm$$
\vspace*{-15pt}$$ \hskip 3cm
+\delta_{i+k,-2\a}\delta_{j+\ell,2\b}[a_2(\a j+\b i)+a'_2(\a+i)]{\bf c}_2,
\eqno(1.3)$$
for $(i,j),(k,\ell)\in Z$.
Note that
${\bf c}_1=0$ if $(-\a,\b)\notin\ZZ$, and
${\bf c}_1={\bf c}_2=0$ if $(-2\a,2\b)\notin\ZZ$.
\par
Set $Z_+=\{(i,j)\in Z\,|\,j\ge-1\}$, $Z_-=\{(i,j)\in Z\,|\,j\le1\}$.
If $\b=-1$, the subspace $\BB^+(\a,-1;a_1,a_2,a'_2)
=\oplus_{(i,j)\in Z_+}\F L_{i,j}\oplus \F{\bf c}_1\oplus \F{\bf c}_2$
is a subalgebra of $\BB(\a,-1;a_1,a_2,a'_2)$.
If $\b=1$, the subspace $\BB^+(\a,1;a_1,a_2,a'_2)
=\oplus_{(i,j)\in Z_-}\F L_{i,j}\oplus \F{\bf c}_1\oplus \F{\bf c}_2$
is a subalgebra of $\BB(\a,1;a_1,a_2,a'_2)$.
\par
For $\a\in\F\bs\{0\}$, we have a Lie algebra $\CC(\a)
=\oplus_{(i,j)\in\sZ\times\sZ}\F L_{i,j}$ with the bracket
$$
[L_{i,j},L_{k,\ell}]=a_{i,j,k,\ell}L_{i+\ell,k+j}
\mbox{ for }i,j,k,\ell\in\Z,
\eqno(1.4)$$
such that
$$
a_{i,j,k,\ell}=\left\{\matrix{
k(j+1)-(\ell+1)i+(\ell-j)\a\hfill&\mbox{if}&j,\ell,j+\ell\ge-1,
\vsp\hfill\cr
[^{\ \,-\ell-2}_{-\ell-j-2}](k(j+1)-(\ell+1)i+(\ell-j)\a)
\hfill&\mbox{if}&j\ge0,\ell\le-2,\vsp\hfill\cr
k-i\hfill&\mbox{if}&j=\ell=-1,\vsp\hfill\cr
-\a+i\hfill&\mbox{if}&j=-1,\ell\le-2,\vsp\hfill\cr
0\hfill&\mbox{if}&j,\ell\le-2,
\hfill\cr}\right.
\eqno(1.5)$$
where $[^i_j]={i!\over j!}$ for $0\le j\le i$ and $[^i_j]=0$ otherwise.
Note that $\CC(\a)$ is a Lie algebra such that $\CC^-(\a)=
\oplus_{i\in\sZ,j\le-2}\F L_{i,j}$ is an abelian ideal and
$\CC(\a)/\CC^-(\a)\cong\BB^+(-\a,-1,1,0,0).$
Certainly we have
the dual Lie algebra $\bar\CC(\a)$ whose $(i,j)$-homogeneous
space is the $(i,-j)$-homogeneous
space of $\CC(\a)$.
\par
The main result of this paper is the following
\par
{\bf Theorem 1.1} {\it
Let $\AA=\oplus_{i,j\in\sZ}\AA_{i,j}$ be a $\ZZ$-graded Lie algebra
 over a field $\F$ of characteristic $0$ with
$\dim\AA_{i,j}\leq 1$ for each $i$ and $j$. Suppose $\AA$ satisfies Conditions (1)-(3).  Then
$\AA$ is isomorphic to Vir$(\a)$
for a suitable $\a\in \F\bs\{0\}$,
or to $\BB(\a,\b;a_1,a_2,a'_2)$ for suitable $\a,\b,a_1,a_2,a'_2\in\F$
with $\a\b\ne0$ and $\b\ne\pm1$, or to
$\BB^+(\a,\pm1;a_1,a_2,a'_2)$ for suitable $\a,a_1,a_2,a'_2\in\F$ with
$\a\ne0$, or to $\CC(\a)$ or to
$\bar\CC(\a)$ for a suitable $\a\in \F\bs\{0\}$.}
\par
For convenience, we combine the two types of algebras Vir$(\a)$
and $\BB(\a,\b)$ into a uniform form: $\EE(\a,\b)=\oplus_{i,j\in Z}
\F L_{i,j}$ subject to
$$
[L_{i,j},L_{k,\ell}]=
(\b (i\ell-jk)+(k-i)+(\ell-j)\a)L_{i+k,j+\ell},
\eqno(1.6)$$
where $\a,\b\in F$ with $\a\ne0$. Note that $\EE(\a,0)=$Vir$(\a)$.
\par\vskip 10pt
\cl{\bf 2. Some basic results on the Virasoro algebra}
\par
Recall that the centerless Virasoro algebra is the Lie algebra
${\rm Vir}=\oz\F L_i$ with bracket $[L_i,L_j]=(j-i)L_{i+j}$ for $i,j\in\Z$.
A weight module $V$ over Vir is called a {\it module of the intermediate
series} if
all the weight multiplicities are $\le1$. There are three families of the
modules of the intermediate series over Vir: $A_{\a,\b},A(\a),B(\a)$ for
$\a,\b\in\F$. They all have basis $\{v_i\,|\,i\in\Z\}$ with the following
actions:
$$
\matrix{
A_{\a,\b}:\!\!\!\!&L_i v_k=(\a+k+\b i)v_{i+k},
\vsp\hfill\cr
A(\a):\!\!\!\!&L_i v_k=(i+k)v_{i+k},\,k\ne0,\
L_i v_0=i(i+\a)v_i,\vsp\hfill\cr
B(\a):\!\!\!\!&L_i v_k=k x_{i+k},\,k\ne-i,\
L_i v_{-i}=-i(i+\a)v_0,\hfill\cr
}
\eqno(2.1)$$
for $i,k\in\Z$. Denote by $A'_{\a,\b}$ the nontrivial irreducible subquotient
module of $A_{\a,\b}$.
\par\vsp
{\bf Theorem 2.1}(see [KS, Z])
{\it (1) A nontrivial irreducible module of the intermediate series over Vir
is isomorphic to $A'_{\a,\b}$ for some suitable
$\a,\b\in\F$.
\vspace*{-3pt}\par
(2) An indecomposable module of the intermediate series over
Vir is isomorphic to a subquotient module of
$A_{\a,\b},A(\a),B(\a)$ for some suitable $\a,\b\in\F$.
\vspace*{-3pt}\par
(3) We have the following isomorphisms:
$$
A_{\a,0}\simeq A_{\a,1}\mbox{ if }\a\notin\Z\mbox{ \ \ and \ \ }
A'_{0,0}\simeq A'_{0,1}.
\eqno(2.2)$$}
\vskip 3pt
\cl{\bf3. Proof of Theorem 1.1}
\par
Now let $\AA=\oplus_{i,j\in\sZ}\AA_{i,j}$ be a $\Z\times\Z$-graded Lie algebra
over $\F$ with $\dim \AA_{i,j}\le 1$
satisfying Conditions (1)-(3). For $i,j\in\Z$, take
$L_{i,j}\in \AA_{i,j}\bs\{0\}$ if $\AA_{i,j}\ne\{0\}$ or else let $L_{i,j}=0$.
By Conditions (1)-(2), we can normalize the elements $L_{i,0}$, $L_{0,\pm1}$
such that
$$
[L_{i,0},L_{j,0}]=(j-i)L_{i+j,0}\hbox{ for }i,j\in\Z.
\eqno(3.1)$$
and
$$
[L_{0,i},L_{0,j}]=(j-i)\a L_{0,i+j}\hbox{ for }-1\le i,j,i+j\le1,
\eqno(3.2)$$
where $\a\in\F\bs\{0\}$ is a constant. Since $\oplus_{i\in\sZ}\AA_{i,\pm1}$ is
an $\AA_0$-module, and $[L_{0,0},L_{0,\pm1}]=\pm\a L_{0,\pm1}\ne0$,
by Theorem 2.3 in [OZ1],
there exist $\b_{\pm1}\in\F$ such that
we can normalize $L_{i,\pm1}$ to satisfy
$$
[L_{i,0},L_{j,\pm1}]=(\pm\a+j+\b_{\pm1}i)L_{i+j,\pm1}
\eqno(3.3)$$
$
\hbox{ for }i,j\!\in\!\Z
\hbox{ with }j\!\pm\!\a,i\!+\!j\!\pm\!\a\!\ne\!0,
 \hbox{ or with } \b_{\pm1}\!\ne\!0,1.
$
In particular $L_{i,\pm1}\ne0$ for $i\in\Z$ with $i\pm\a\ne0$.
\par
Now we divide the proof of Theorem 1.1 into several lemmas.
\par
{\bf Lemma 0.} After a change of basis, the coefficients $\b_1$ and $\b_{-1}$
defined in (3.3) satisfy the relation
$$
\b_1\in\{\b_{-1},-\b_{-1},-1-\b_{-1},-2-\b_{-1}\}.
$$
\par
{\bf Proof.}
Suppose
$$
[L_{i,-1},L_{j,1}]=c_{i,j}L_{i+j,0}\hbox{ for }i,j\in\Z
\hbox{ and some }c_{i,j}\in\F,
\eqno(3.4)$$
where $c_{0,0}=2\a.$ Applying ${\rm ad\ssc\,}L_{k,0}$ to (3.4), using
(3.3), we have that
$$
(-\a+i+\b_{-1}k)c_{i+k,j}+(\a+j+\b_1k)c_{i,j+k}=(i+j-k)c_{i,j},
\eqno(3.5)$$
holds for all $i,j,k\in\Z$ with
$$
(i-\a)(i+k-\a)\ne0\hbox{ or }\b_{-1}\ne0,1,\hbox{ and }
(j+\a)(j+k+\a)\ne0\hbox{ or }\b_1\ne0,1.
\eqno(3.6)$$
It is clear that $\b_1$ and $\b_{-1}$ are
symmetric in (3.5). We shall solve $c_{i,j}$ in (3.5).
Setting $i=0,j=k$ and $j=0,i=k$ in (3.5), we obtain respectively
$$
\matrix{
(-\a+\b_{-1}k)c_{k,k}+(\a+(1+\b_1)k)c_{0,2k}=0
\hbox{ if }(k-\a)(k+\a)(2k+\a)\ne0,
\vsp\hfill\cr
(-\a+(1+\b_{-1})k)c_{2k,0}+(\a+\b_1k)c_{k,k}=0
\hbox{ if }(k-\a)(2k-\a)(k+\a)\ne0.
\hfill\cr}
\eqno(3.7)$$
{}From (3.7), we obtain
$$
(-\a+\b_{-1}k)(-\a+(1+\b_{-1})k)c_{2k,0}-
(\a+\b_1k)(\a+(1+\b_1)k)c_{0,2k}=0,
\eqno(3.8)$$
for $k\in\Z$ with $k\pm\a,2k\pm\a\ne0$.
In (3.5), setting $i=0$ and substituting $j,k$ by $2j,2i$, and
setting $j=0$ and substituting $i,k$ by $2i,2j$ respectively, we obtain
$$
\matrix{
(-\a+2\b_{-1}i)c_{2i,2j}+(\a+2j+2\b_1i)c_{0,2(i+j)}=2(j-i)c_{0,2j},
\vsp\hfill\cr
(-\a+2i+2\b_{-1}j)c_{2(i+j),0}+(\a+2\b_1j)c_{2i,2j}=2(i-j)c_{2i,0},
\hfill\cr}
\eqno(3.9)$$
for $2i-\a,2j+\a,2(i+j)\pm\a\ne0$. From (3.9), we obtain
$$
\matrix{
(\a+2\b_1j)((\a+2j+2\b_1i)c_{0,2(i+j)}-2(j-i)c_{0,2j})
\vsp\hfill\cr
=(-\a+2\b_{-1}i)((-\a+2i+2\b_{-1}j)c_{2(i+j),0}-2(i-j)c_{2i,0}).
\hfill\cr}
\eqno(3.10)$$
In (3.5), setting $i=j=0$ and substituting $k$ by $2k$, we obtain
$$
(-\a+2\b_{-1}k)c_{2k,0}+(\a+2\b_1k)c_{0,2k}=-4k\a
\hbox{ for }2k\pm\a\ne0.
\eqno(3.11)$$
Multiplying (3.10) by $-\a+2\b_{-1}(i+j)$ and using (3.11) to substitute
$c_{2(i+j),0},c_{2i,0}$, we obtain
$$
\matrix{
(-\a+2\b_{-1}(i+j))(\a+2\b_1j)((\a+2j+2\b_1i)c_{0,2(i+j)}-2(j-i)c_{0,2j})
\vsp\hfill\cr
=(-\a+2\b_{-1}i)(-\a+2i+2\b_{-1}j)(-(\a+2\b_1(i+j))c_{0,2(i+j)}-4(i+j)\a)
\vsp\hfill\cr
-2(j-i)(-\a+2\b_{-1}(i+j))((\a+2\b_1i)c_{0,2i}+4i\a),
\hfill\cr}
\eqno(3.12)$$
for $2i\pm\a,2j\pm\a,2(i+j)\pm\a\ne0$.
Setting $i=-j$ gives
$$
\matrix{
-\a(\a-2\b_1i)((\a+2(\b_1-1)i)2\a+4ic_{0,-2i})
\vsp\hfill\cr
=-\a(-\a+2\b_{-1}i)(-\a+2(1-\b_{-1})i)2\a
-4i\a((\a+2\b_1i)c_{0,2i}+4i\a),
\hfill\cr}
\eqno(3.13)$$
for $\pm 2i\pm\a\ne0$, that is
$$
(\a+2\b_1i)c_{0,2i}-(\a-2\b_1i)c_{0,-2i}
=2i\a((\b_{-1}-\b_1)(\b_{-1}+\b_1-1)-2),
\eqno(3.14)$$
for $i,\pm2i\pm\a\ne0$.
{}From (3.8), (3.11), we have
$$
\matrix{
\left|\matrix{
(-\a\!+\!\b_{-1}k)(-\a\!+\!(1\!+\!\b_{-1})k)&\!\!\!
-(\a\!+\!\b_1k)(\a\!+\!(1\!+\!\b_1)k)\vsp\cr
-\a+2\b_{-1}k&\!\!\!\a+2\b_1k\cr}\right|
c_{0,2k}
\vspace*{6pt}\hfill\cr\hskip 8cm
\!=\!-4k(-\a\!+\!\b_{-1}k)(-\a\!+\!(1\!+\!\b_{-1})k)\a,
\hfill\cr}
\eqno(3.15)$$
for $k,k\pm\a,2k\pm\a\ne0.$ Denote by $d_k$ the determinant in (3.15).
Using (3.15) in (3.14), we obtain
$$
\matrix{
2(-\a+\b_{-1}i)(-\a+(1+\b_{-1})i)
(\a+2\b_1i)d_{-i}
\vsp\hfill\cr\ \ \ \
+2(\a\!+\!\b_{-1}i)(\a\!+\!(1\!+\!\b_{-1})i)(\a\!-\!2\b_1i)d_i
\!+\!((\b_{-1}\!-\!\b_1)(\b_{-1}\!+\!\b_1\!-\!1)\!-\!2)d_id_{-i}\!=\!0,
\hfill\cr}
\eqno(3.16)$$
for $\pm i,\pm i\pm\a,\pm2i\pm\a\ne0.$
First assume that one of $\b_1,\b_{-1}$ is zero, say,
$\b_{-1}=0$ (since $\b_1$ and $\b_{-1}$ are symmetric).
Then from (3.15), we see that $d_i$ is a polynomial on $i$ of degree 2 with
the coefficient of $i^2$ being $-\a\b_1(\b_1+3)$,
and (3.16) is a polynomial in $i$ of degree 4.
Calculating the coefficient of $i^4$ gives
$$
\b^2_1(\b_1+1)(\b_1+2)(\b_1-1)(\b_1+3)\a^2=0\hbox{ if }\b_{-1}=0.
\eqno(3.17)$$
Now assume that $\b_1\b_{-1}\ne0$.
Then (3.16) is a polynomial on $i$ of degree 6.
Calculating the coefficient of $i^6$ gives
$$
-4(\b_1-\b_{-1})(\b_1+\b_{-1})(\b_1+\b_{-1}+1)(\b_1+\b_{-1}+2)
\b_1^2\b_{-1}^2=0
\hbox{ if }\b_1,\b_{-1}\ne0.
\eqno(3.18)$$
{}From (3.17), (3.18), we obtain that
$$
\b_1\!=\!\b_{-1},-\b_{-1},-1\!-\!\b_{-1},-2\!-\!\b_{-1},
\hbox{ or }\b_1\!=\!1,-3,\,\b_{-1}\!=\!0,\hbox{ or }\b_1\!=\!0,\,
\b_{-1}\!=\!1,-3.
\eqno(3.19)$$
If $\b_1=1,-3$ and $\b_{-1}=0$, we can re-choose
$\{L_{i,-1}\,|\,i\in\Z\bs\{\a\}\}$ (if $\a\in\Z$, we shall not consider
$L_{\a,-1}$) such that $\b_{-1}$ becomes 1 (cf.~(2.2)). By interchanging
$\AA_{0,1}$ with $\AA_{0,-1}$ if necessary, it suffices to consider the
following 4 cases:
$\b_1\!=\!\b_{-1},-\b_{-1},-1\!-\!\b_{-1},-2\!-\!\b_{-1}$.
\qed\par
{\bf Remark.}
Before we consider cases one by one, observe the following.
{}From (3.15), we can solve $c_{0,2k}$, then from the first equation of
(3.9), we can solve $c_{2i,2j}$. Note that by setting $k=\pm1$ in
(3.5), we can express $c_{i,j}$ in terms of $c_{i-1,j+1},c_{i+1,j-1}$.
So if $i+j$ is even, we can solve $c_{i,j}$. On the other hand, if
$i+j$ is odd, by taking $k$ to be odd in (3.5), we can solve $c_{i,j}$.
This shows that the solution of $c_{i,j}$ in (3.5) is unique (except
possibly a finite number of $i,j$).
\par
{\bf Lemma 1.} If $\b_1=\b_{-1}$, then $\b_1=0,\pm1,-{1\over2}$. The cases
$\b_1=0,1,-{1\over2}$ are covered by the cases
$\b_1\in\{-\b_{-1},-1-\b_{-1},-2-\b_{-1}\}$ considered below, while the case
$\b_1=\b_{-1}=1$ leads to a contradiction.
\par
{\bf Proof.}
In this case, (3.15) and (3.11) give
$$
c_{0,2k}=
{2\a(-\a+\b_1 k)(-\a+(1+\b_1)k)\over
\a^2-2\b_1^2(1+\b_1)k^2},\ \
c_{2k,0}=
{2\a(\a+\b_1 k)(\a+(1+\b_1)k)\over
\a^2-2\b_1^2(1+\b_1)k^2},
\eqno(3.20)$$
for $k,k\pm\a,2k\pm\a\ne0$ and $\a^2-2\b_1^2(1+\b_1)k^2\ne0$.
Using this in (3.10), we obtain
$$
{16\a^2\b_1^3(\b_1-1)(1+\b_1)(1+2\b_1)ij(i+j)(j-i)^2(\a^2+\b_1^2(1+\b_1)ij)
\over
(\a^2-2\b_1^2(1+\b_1)i^2)(\a^2-2\b_1^2(1+\b_1)j^2)
(\a^2-2\b_1^2(1+\b_1)(i+j)^2)}=0,
\eqno(3.21)$$
for all but a finite number of $i,j\in\Z$.
Thus $\b_1=0,\pm1,-{1\over2}$.
Cases $\b_1=0,-{1\over2},-1$ will be covered by cases
$\b_1=-\b_{-1},-1-\b_{-1},-2-\b_{-1}$ considered below.
Thus suppose $\b_1=1$.
By (3.20) and the first equation of (3.9), we obtain
$$
c_{2i,2j}={2\a(\a^2+\a(i-j)-4ij)\over(\a+2j)(\a-2i)},
\eqno(3.22)$$
for all but a finite number of $i,j\in\Z$.
Setting $i,j,k$ in (3.5) to be $2i,2j,2k$ respectively,
using (3.22), we see that
(3.5) does not hold. Thus we obtain the lemma.
\qed\par
{\bf Lemma 2.} If $\b_1=-\b_{-1}$, then $\a\notin\Z$, and after a change of
notation, we have
$\b_1=0$ which is covered by the cases $\b_1=-1-\b_{-1}$ considered
below.
\par
{\bf Proof.}
By (3.15), (3.11) and the first equation of (3.9), we have
$$
\matrix{\dis
c_{0,2k}={2\a(\a+(\b_1-1)k)\over\a+2\b_1k},\ \
c_{2k,0}={2\a(\a+(\b_1+1)k)\over\a+2\b_1k},
\vspace*{4pt}\hfill\cr\dis
c_{2k,2k}=
{2\a(\a+2(\b_1-1)k)(\a+2(\b_1+1)k)\over
(\a+2\b_1k)(\a+4\b_1k)},
\hfill\cr}
\eqno(3.23)$$
for all but a finite number of $k\in\Z$. Setting $i,j,k$ in (3.5)
to be $2k,2k,-2k$ respectively, using (3.23), we obtain
$$
-{48\a(\b_1-1)(\b_1+1)k^3\over(\a+2\b_1k)(\a+4\b_1k)}=0,
\eqno(3.24)$$
for all but a finite number of $k\in\Z$.
This proves $\b_1=\pm1$. By symmetry, we can suppose $\b_1=1,\b_{-1}=-1$.
Assume that $\a\in\Z$. Choose $0\ne k\in\Z$ such that $k-\a$ is even.
Setting $j=-\a,i=\a+k$ in (3.5), we obtain that $c_{\a+k,k-\a}=0$
for all $k\ne0$. In particular, $c_{2\a,0}=0$. But from (3.23) we see that
$c_{2\a,0}=2\a\ne0$, which is
 a contradiction. Thus
$\a\notin\Z$. Thus by (2.2), we can re-take $\b_1$ to be zero.
\qed\par
{\bf Lemma 3.} If $\b_1=-1-\b_{-1}$, then after a change of notation, we
have
$\b_1=-2,{1\over2}$. The case $\b_1=-2$ is covered by the case
$\b_1=-2-\b_{-1}$ considered below,
while the case $\b_1={1\over2}$ leads to a
contradiction.
\par
{\bf Proof.}
By interchanging $\AA_{0,1}$ with $\AA_{0,-1}$ if necessary, we can suppose
$\b_1\ne 0,1$ since $\b_1+\b_{-1}=-1$. So $L_{i,1}\ne0$ for all $i\in\Z$.
{}From (3.15) we get $c_{0,2i}=2\a$ for almost all $i\in \Z$. Then from (3.9)
and (3.5), we get
$$
c_{i,j}=2\a\mbox{ for }i,j\in\Z
\mbox{ (except possibly one $i$ if $\a\in\Z,\,\b_{-1}=0,1$)}.
\eqno(3.25)$$
For $i,j,k\in\Z$, we have
$$
[L_{k,-1},[L_{i,1},L_{j,1}]]
=c_{k,i}[L_{k+i,0},L_{j,1}]-c_{k,j}[L_{k+j,0},L_{i,1}]
=2\a(j-i)(1-\b_1)L_{i+j+k,1}.
\eqno(3.26)$$
This shows that $[L_{i,1},L_{j,1}]\ne0$ for $i\ne j$. Thus $\AA_{j,2}
\ne\{0\}$ for all $j\in\Z$. Choose a basis $L_{j,2}$ of $\AA_{j,2}$ as
follows:
$$
L_{j,2}=\left\{\matrix{\dis{1\over j}[L_{0,1},L_{j,1}]\hfill\!\!&
\mbox{if }j\ne0,\vsp\hfill\cr\dis{1\over2}[L_{-1,1},L_{1,1}]
\hfill\!\!&\mbox{if }j=0,\hfill\cr}\right.
\eqno(3.27)$$
Then by (3.26),
$$
[L_{k,-1},L_{i,2}]=2\a(1-\b_1)L_{i+k,1}\mbox{ for }i,k\in\Z.
\eqno(3.28)$$
Since $[L_{i,1},L_{j,1}]\in\AA_{i+j,2}$, assume that $[L_{i,1},L_{j,1}]=
c'_{i,j}L_{i+j,2}$ for some $c'_{i,j}\in\F$,
then using (3.26), (3.28), we obtain
$$
[L_{i,1},L_{j,1}]=(j-i)L_{i+j,2}\mbox{ for }i,j\in\Z.
\eqno(3.29)$$
Applying ${\rm ad\ssc\,}L_{k,0}$ to (3.29), we obtain
$$
[L_{k,0},L_{i,2}]=(2\a+i+(2\b_1+1)k)L_{i+k,2}\mbox{ for }i,k\in\Z.
\eqno(3.30)$$
Now by (3.25), (3.28)-(3.30), we have
$$
\matrix{
[L_{k,-1},[L_{i,1},L_{j,2}]]
\!\!\!\!&
=2\a([L_{i+k,0},L_{j,2}]+(1-\b_1)[L_{i,1},L_{j+k,1}])
\vsp\hfill\cr&
=2\a(2\a+3\b_1 i+(2-\b_1)j+(2+\b_1)k)L_{i+j+k,2}\mbox{ for }i,j,k\in\Z.
\hfill\cr}
\eqno(3.31)$$
Thus $\AA_{i,3}\ne\{0\}$ for all $i\in\Z$. Since $\oplus_{i\in\sZ}\AA_{i,3}$
is a Vir-module, we can choose $0\ne L_{i,3}\in\AA_{i,3}$ for
$i\in\Z$ such that
$$
[L_{k,0},L_{i,3}]=(3\a+i+\b_3 k)L_{i+k,3}
\mbox{ for }i,k\!\in\!\Z\mbox{ (and $i\!+\!\a,
i\!+\!k\!+\!\a\!\ne\!0$ if $\b_3\!=\!0,1$)},
\eqno(3.32)$$
where $3\a$ appears because it is the eigenvalue of $L_{0,0}$ on
$L_{i,3}$ (cf.~(2.1)).
Suppose
$$
[L_{i,1},L_{j,2}]=d_{i,j}L_{i+j,3}\mbox{ for }i,j\in\Z
\mbox{ and some }d_{i,j}\in\F,
\eqno(3.33)$$
such that at least some $d_{i,j}\ne0$ by (3.31).
Applying ${\rm ad\ssc\,}L_{k,-1}$ to (3.33), by (3.31), we obtain
$$
d_{i,j}[L_{k,-1},L_{i+j,3}]=2\a(2\a+3\b_1 i+(2-\b_1)j+(2+\b_1)k)L_{i+j+k,2}
\mbox{ for }i,j,k\in\Z.
\eqno(3.34)$$
Setting $i=0$ and replacing $j$ by $i+j$ in (3.34), we obtain
$$
d_{0,i+j}[L_{k,-1},L_{i+j,3}]=2\a(2\a+(2-\b_1)(i+j)+(2+\b_1)k)L_{i+j+k,2}
\mbox{ for }i,j,k\in\Z.
\eqno(3.35)$$
Formulas (3.34), (3.35) show that
$$
d_{i,j}(2\a\!+\!(2\!-\!\b_1)(i\!+\!j)\!+\!(2\!+\!\b_1)k)
\!=\!d_{0,i+j}(2\a\!+\!3\b_1 i
\!+\!(2\!-\!\b_1)j\!+\!(2\!+\!\b_1)k)
\mbox{ for }i,j,k\!\in\!\Z.
\eqno(3.36)$$
By (3.3), (3.30), (3.32), (3.33), we have
$$
\matrix{
d_{0,i}(3\a+i+\b_3k)L_{3,i+k}
\!\!\!\!&=[L_{k,0},[L_{0,1},L_{i,2}]]
\vspace*{4pt}\hfill\cr&
=((\a+\b_1k)d_{k,i}+(2\a+i+(2\b_1+1)k)d_{0,i+k})
L_{3,i+k}.\hfill\cr}
\eqno(3.37)$$
Comparing the coefficient of $k$ in (3.36)
gives $d_{i,j}=d_{0,i+j}$ or $2+\b_1=0$.
Assume that $\b_1=-2$. Then $\b_{-1}=1$, and by (2.2),
we can re-choose $\b_{-1}=0$. Thus this becomes a special case of Case 4
considered below.
Thus we can suppose $d_{i,j}=d_{0,i+j}$. Then (3.36) gives $\b_1={1\over2}$.
Then (3.37) gives
$$
d_{0,i+k}=
{d_{0,i}(3\a+i+\b_3k)\over
3\a+i+{5k\over2}}.
\eqno(3.38)$$
In particular
$$
d_{0,k}=
{d_{0,0}(3\a+\b_3k)\over
3\a+{5k\over2}}.
\eqno(3.39)$$
Substituting (3.39) into (3.38) we obtain that
$$
(3\a+\b_3 i)(3\a+i+\b_3 k)(3\a+{5(k+i)\over2})=
(3\a+\b_3(k+i))(3\a+i+{5k\over2})(3\a+{5i\over2}),
\eqno(3.40)$$
i.e.,
$${1\over4} (-5 + 2 \b_3) i k (5 \b_3 k + 5 \b_3 i + 6 \b_3 \a + 9 \a)=0.$$
Then $\b_3={5\over2}$ and $d_{0,i}=d_{0,0}$ is a constant.
By re-choosing $L_{i,3}$ for all $i\in\Z$,
 we can assume $d_{i,j}=1$, and (3.33) means that the left-hand side
of (3.31) is $[L_{k,-1},L_{i+j,3}]$.
Note that we have got
$$\b_1={1\over2},\,\,\b_{-1}=-{3\over2},\,\,\b_3={5\over2}.
\eqno(3.41)$$
\par
Similar to (3.26), we have
$$
[L_{k,1},[L_{i,-1},L_{j,-1}]]=-5\a(j-i)L_{i+j+k,-1}\mbox{ for }i,j,k\in\Z,
$$
which shows that $[L_{i,-1},L_{j,-1}]\ne0$ for $i\ne j$. Thus $\AA_{j,-2}
\ne\{0\}$ for all $j\in\Z$. Choose a basis $L_{j,-2}$ of $\AA_{j,-2}$
to be
$$
L_{j,-2}={1\over j}[L_{j,-1},L_{0,-1}]\mbox{ if }j\ne0\mbox{ and }
L_{0,-2}={1\over2}[L_{1,-1},L_{-1,-1}],
\eqno(3.42)$$
then
$$
[L_{i,-1},L_{j,-1}]=(i-j)L_{i+j,-2},
[L_{i,1},L_{j,-2}]=5\a L_{i+j,-1},
[L_{i,0},L_{j,-2}]=(-2\a+j-2i)L_{i+j,-2},
$$
and
$$
[L_{k,1},[L_{i,-1},L_{j,-2}]]=\a(4\a+9i-7j-k)L_{i+j+k,-2}.
\eqno(3.43)$$
This shows that $\AA_{i,-3}\ne\{0\}$ for $i\in\Z$. Choose
$L_{i,-3}\in\AA_{i,-3}\bs\{0\}$ and suppose
$[L_{i,-1},L_{j,-2}]=d'_{i,j}L_{i+j,-3}$ for some $d'_{i,j}\in\F$.
Then at least some $d'_{i,j}\ne0$.
Setting $i=0$ and replacing $j$ by $i+j$ in (3.43), comparing the result
with (3.43), we obtain
$$
d'_{i,j}(4\a-7i-7j-k)=d'_{0,i+j}(4\a+9i-7j-k)\mbox{ for }i,j,k\in\Z.
\eqno(3.44)$$
Comparing the coefficient of $k$ we get $d'_{i,j}=d'_{0,i+j}$, and then we
obtain that (3.44) cannot hold for all $i,j,k\in\Z$.
Thus $\b_1={1\over2}$ cannot occur.
\qed\par
Therefore, it remains to consider the only case $\b_1=-2-\b_{-1}$.
By interchanging $\AA_{0,1}$ with $\AA_{0,-1}$ if necessary, we can suppose
$\b_1\ne 0$.
We denote $\b=\b_1+1$. Then $\b\ne1$.
\par
{\bf Lemma 4.1.} If $\b_1=-2-\b_{-1}$ and $\b_{-1}\ne0$, then
$\AA=\EE(\a,\b)$ or $\BB(\a,\b;a_1,a_2,a'_2)$.
\par
{\bf Proof.}
First we claim:
one can choose $\{L_{i,j}\,|\,(i,j)\in Z^+\}$ such that
(1.6) holds for $(i,j),(k,\ell),(i+k,j+\ell)\in Z^+$.
\par
Using (3.15) and the statements after (3.19), we get
$$
c_{i,j}=2\a+(\b-1)i+(\b+1)j
\mbox{ for }(i,1),(j,1)\in Z.
\eqno(3.45)$$
For $(j,2)\in Z$, choose $i,k\in\Z$ such that
$(k,-1),(i,1),(j-i,1),(j+k,2)\in Z$, then
$$
[L_{k,-1},[L_{i,1},L_{j-i,1}]]
=(1-\b)((2\b-1)k+3\a+(\b+1)j)(j-2i)L_{j+k,1}.
\eqno(3.46)$$
This shows that if $j\ne2i$, then $0\ne[L_{i,1},L_{j-i,1}]\in\AA_{j,2}$.
So $\AA_{j,2}\ne0$ for all $(j,2)\in Z$. For $(j,2)\in Z$,
choose minimal $i\in\Z_+$ such that $2i\ne j$ and $(i,1),(j-i,1)\in Z$ and
define
$$
L_{j,2}={1\over (1-\b)(j-2i)}[L_{i,1},L_{j-i,1}].
\eqno(3.47)$$
Then by (3.46),
$$
[L_{k,-1},L_{i,2}]=
((2\b-1)k+3\a+(\b+1)i)L_{i+k,1}\mbox{ for }
(k,-1),(i,2),(i+k,1)\in Z.
\eqno(3.48)$$
Assume that $[L_{i,1},L_{j,1}]=c'_{i,j}L_{i+j,2}$ for
$(i,1),(j,1),(i+j,2)\in Z$,
then by (3.46)-(3.48), $c'_{i,j}=(1-\b)(j-i)$, i.e.,
$$
[L_{i,1},L_{j,1}]=(1-\b)(j-i)L_{i+j,2}\mbox{ for }(i,1),(j,1),(i+j,2)\in Z.
\eqno(3.49)$$
Applying ${\rm ad\ssc\,}L_{k,0}$ to it, using
(3.3), (3.47), we can deduce
$$
[L_{k,0},L_{j,2}]=(2\a+j+(2\b-1)k)L_{j+k,2}
\mbox{ for }k\in\Z,(j,2),(j+k,2)\in Z.
\eqno(3.50)$$
For $m\ge2$, assume that we have chosen $L_{i,m}$ for $(i,m)\in Z$ such
that (1.6) holds with $-1\le j,\ell,j+\ell\le m$.
For $(j,m+1)\in Z$, choose $i,k\in\Z$ such that
$$
(k,-1),(i,1),(j-i,m),(j-i+k,m-1),(j+k,m)\in Z,
\eqno(3.51)$$
then
$$
\matrix{
[L_{k,-1},[L_{i,1},L_{j-i,m}]]
\vsp\hfill\cr
=\!(\a(m\!+\!2)\!+\!k(\b(m\!+\!1)\!-\!1)
\!+\!(\b\!+\!1)j)(\a(m\!-\!1)\!+\!(1\!-\!\b)j\!+\!(\b(m\!+\!1)\!-\!2)i)
L_{j+k,m}.
\hfill\cr}
\eqno(3.52)$$
This shows that $0\ne[L_{i,1},L_{j-i,m}]\in\AA_{j,m+1}$.
So $\AA_{j,m+1}\ne0$ for all $(j,m+1)\in Z$. For $(j,m+1)\in Z$,
choose minimal $i\in\Z_+$ such that (3.51) holds with $k=0$, then
we define
$$
L_{j,m+1}=
(\a(m-1)+(1-\b)j+(\b(m+1)-2)i)^{-1}[L_{i,1},L_{j-i,m}].
\eqno(3.53)$$
Then as in the arguments after (3.47), we obtain that
(1.6) holds with $-1\le j,\ell,j+\ell\le m+1$.
This proves the claim.
\par
Similarly, as in the proof of the claim above, we have the following claim:
one can choose $\{L_{i,j}\,|\,(i,j)\in Z\}$ such that (1.6) holds.
\par
So if $Z=\ZZ$ or $\AA_{i,j}=0$ for $(i,j)\notin Z$,
then we obtained that $\AA=\EE(\a,\b)$.
Thus assume that $Z\ne\ZZ$, and $\AA_{i_0,j_0}\ne0$ for some
$(i_0,j_0)\in\ZZ\bs Z$.
Since $\AA$ is generated by $\AA_0$ and
$\AA_{0,\pm1}$, we can find some $i,j\in\Z$ such that
$(i,j),(i_0-i,j_0-j)\in Z$ and
$$
L_{i_0,j_0}=[L_{i,j},L_{i_0-i,j_0-j}]\ne0.
\eqno(3.54)$$
Then by applying ${\rm ad\,}L_{i,0}$ and
${\rm ad\,}L_{0,\pm1}$ to (3.54), we see
that $L_{i_0,j_0}$ is a central element of $\AA$. Since central extensions
are determined by 2-cocycles, by [DZ] we obtain that $\AA$ has
the form (1.3).
\qed\par
{\bf Lemma 4.2a.} If $\b_1=-2-\b_{-1}$ and $\b_{-1}=0$ with
${\rm dim}(\oplus_{i\in\sZ}\AA_{i,-2})\ge2$, then $\AA=\CC(\a)$.
\par
{\bf Proof.}
Then $\b_1=-2$ and $\b=-1$. By (3.3), (3.46), we have
$$
[L_{k,1},[L_{i,-1},L_{j,-1}]]=0\mbox{ for }i,j,k\in\Z.
\eqno(3.55)$$
First assume that ${\rm dim}(\oplus_{i\in\sZ}\AA_{i,-2})\ge2$.
Then $\oplus_{i\in\sZ}\AA_{i,-2}$ is a nontrivial $\Vir$-module.
Thus ${\rm dim}(\oplus_{i\in\sZ}\AA_{i,-2})=\infty$ and
there exists $\b_{-2}\in\F$ such that
we can choose suitable $L_{i,-2}\in\AA_{i,-2}\bs\{0\}$
if $\AA_{i,-2}\ne\{0\}$ and $L_{i,-2}=0$ if $\AA_{i,-2}=\{0\}$
such that
$$
[L_{k,0},L_{i,-2}]=(-2\a+i+\b_{-2}k)L_{i+k,-2},
\eqno(3.56)$$
for all $k\in\Z$ and all but a finite number of $i\in\Z$
(note that (3.56) holds for all $i,k\in\Z$ if $2\a\notin\Z$ or $\b_{-2}\ne0,1$,
and (3.56) holds for $i,i+k\ne2\a$ if $2\a\in\Z$ and $\b_{-2}=0,1$).
Since $\oplus_{i\in\sZ}\AA_{i,-2}=
\oplus_{i,j\in\sZ}[\AA_{i,-1},\AA_{j,-1}]$, by (3.55) we have
$$
[L_{i,1},L_{j,-2}]=0\mbox{ for }i,j\in\Z.
\eqno(3.57)$$
{}From this and (3.45), (3.56), we have
$$
[L_{k,1},[L_{i,-1},L_{j,-2}]]=
-(2\a-2i)(-2\a+j+(i+k)\b_{-2})L_{i+j+k,-2}.
\eqno(3.58)$$
Thus $\AA_{i,-3}\ne\{0\}$ for all $i\in\Z$.
Thus there exist $\b_{-3}\in\F$ and $L_{i,-3}\in\AA_{i,-3}\bs\{0\}$
such that
$$
[L_{k,0},L_{i,-3}]=(-3\a+i+\b_{-3}k)L_{i+k,-3},
\eqno(3.59)$$
for all $i,k\in\Z$ with $(i-3\a)(i+k-3\a)\ne0$.
Write
$$
[L_{i,-1},L_{j,-2}]=d_{i,j}L_{i+j,-3}\mbox{ for }i,j\in\Z
\mbox{ and some }d_{i,j}\in\F,
\eqno(3.60)$$
then at least some $d_{i,j}\ne0$.
Applying ${\rm ad\ssc\,}L_{k,0}$ to (3.60), we obtain
$$
(-\a+i)d_{i+k,j}+(-2\a+j+k)d_{i,j+k}=(-3\a+i+j+\b_{-3}k)d_{i,j},
\eqno(3.61)$$
for $i,j,k\in\Z$ with $(j+k-3\a)(j-3\a)(i+j-3\a)(i+j+k-3\a)\ne0$.
Setting $i=0$ and replacing $j$ by $i+j$ in (3.58), comparing the result
with (3.58), we obtain
$$
2\a(-2\a+i+j+\b_{-2}k)d_{i,j}=(2\a-2i)(-2\a+j+(i+k)\b_{-2})d_{0,i+j},
\eqno(3.62)$$
for all but a finite number of $i,j,k\in\Z$.
This shows that $\b_{-2}=0$ or $d_{i,j}=\a^{-1}(\a-i)d_{0,i+j}$.
Using (2.2), we can suppose $\b_{-2}=1$.
Then we have $d_{i,j}=\a^{-1}(\a-i)d_{0,i+j}$.
Using this in (3.61),
we get that $\b_{-3}=2$ and $d_{0,i}$ is a constant. Thus we
can take $d_{0,i}
=-\a$ to get $d_{i,j}=-\a+i$.
Write
$$
[L_{i,-1},L_{j,-1}]=d'_{i,j}L_{i+j,-2}\mbox{ for }i,j\in\Z
\mbox{ and some }d'_{i,j}\in\F,
\eqno(3.63)$$
such that some $d'_{i,j}\ne0$.
Applying ${\rm ad\ssc\,}L_{k,0}$ to this equation, using (3.56),
we obtain
$$
(-\a+i)d'_{i+k,j}+(-\a+j)d'_{i,j+k}=(-2\a+i+j+k)d'_{i,j},
\eqno(3.64)$$
for all but a finite number of $i,j,k\in\Z$.
Note that $d'_{i,i}=0$ and it is easy to
see that this system has a unique solution
up to a nonzero scalar and $d'_{i,j}=j-i$ is a solution. Thus we can take
$d'_{i,j}=j-i$ for $i,j\in\Z$. If $\a\in\Z$, by applying
${\rm ad\ssc\,}L_{\a-i,0}$ to (3.63), we see that $L_{\a,-1}\ne0$ and
$[L_{\a-i,0},L_{i,-1}]=(-\a+i)L_{\a,-1}$ for $i\in\Z$. Thus, as
a $\Vir$-module, $\oplus_{i\in\sZ}\AA_{i,-1}$ is isomorphic to $A_{-\a,0}$
(cf.~(2.1)).
Similarly if $2\a\in\Z$, by considering
$[L_{k,-1},[L_{i,-1},L_{2\a-i,-1}]]$ and using (3.60),
we can obtain $L_{2\a,-2}\ne0$ and
$\oplus_{i\in\sZ}\AA_{i,-2}\cong A_{-2\a,1}$ as $\Vir$-modules.
Thus we have proved (1.6) for $j\in\{0,\pm1\}$, $\ell\ge -2$
and for all $i,k\in\Z$.
\par
Our first claim is: Formula (1.6) holds for $j\in\{0,\pm1\}$
and for all $i,k,\ell\in\Z$.
\par
We shall show this claim by induction on $-\ell$.  We know that
the claim holds for $\ell\ge-2$. Assume that the claim holds for
$\ell\le-2$.
Then we have
$$
\matrix{
[L_{k,1},\!\!\!\!\!&[L_{i,-1},L_{j,\ell}]]
\vsp\hfill\cr&
=
((-2\a\!+\!2i)[L_{i+k,0},L_{j,\ell}]\!-\!(\ell\!+\!2)((\ell\!-\!1)\a\!+\!2j\!-
\!(\ell\!+\!1)k)[L_{i,-1},L_{j+k,\ell+1}]
\vsp\hfill\cr&
=(\a-i)(\ell+1)((\ell-2)\a+2(i+j)-\ell k)L_{i+j+k,\ell}.
\hfill\cr}
\eqno(3.65)$$
Thus $\AA_{i,\ell-1}\ne\{0\}$ for all $i\in\Z$ and so
there exist $\b_{\ell-1}\in\F$ and $L_{i,\ell-1}\in\AA_{i,\ell-1}\bs\{0\}$
such that
$$
[L_{k,0},L_{i,\ell-1}]=((\ell-1)\a+i+\b_{\ell-1}k)L_{i+k,\ell-1},
\eqno(3.66)$$
for all $k\in\Z$ and all but a finite number of $i\in\Z$.
Write
$$
[L_{i,-1},L_{j,\ell}]=d^{(\ell)}_{i,j}L_{i+j,\ell-1}\mbox{ for }i,j\in\Z
\mbox{ and some }d^{(\ell)}_{i,j}\in\F,
\eqno(3.67)$$
then at least some $d^{(\ell)}_{i,j}\ne0$.
Applying ${\rm ad\ssc\,}L_{k,0}$ to (3.67), we obtain
$$
(-\a+i)d^{(\ell)}_{i+k,j}+(\ell\a+j-(\ell+1)k)d^{(\ell)}_{i,j+k}=
((\ell-1)\a+i+j+\b_{\ell-1}k)d^{(\ell)}_{i,j},
\eqno(3.68)$$
for all but a finite number of $i,j,k\in\Z$.
Setting $i=0$ and replacing $j$ by $i+j$ in (3.65), comparing the result
with (3.65), we obtain $d_{i,j}^{(\ell)}=\a^{-1}(\a-i)d_{0,i+j}^{(\ell)}$
for $i,j\in\Z$.
Using this in (3.68), we obtain
$\b_{\ell-1}=-\ell$ (thus (3.66), (3.68) hold for all $i,j,k\in\Z$)
and $d^{(\ell)}_{0,i}$ is a constant. Thus we can set
$d^{(\ell)}_{0,i}=-\a$ to give $d^{(\ell)}_{i,j}=-\a+i$.
Then (3.65) and (3.67) show that
$$
[L_{i,1},L_{k,\ell-1}]=
(\ell+1)((\ell-2)\a+2k-\ell i)L_{i+k,\ell}
\mbox{ for }i,k,\ell\in\Z.
\eqno(3.69)$$
This proves the claim.
\par
Our second claim is: Formula (1.5) holds for  all $i,j,k,\ell\in\Z$.
\par
Formula (1.5) holds if $j,\ell\ge-1$, or $j=0,\pm1$ and
$\ell\in\Z$. By skew-symmetry, we only need to show
(1.5) for $j\le-2$. So we assume that $j\le-2$.
\par
If $j=-2$, then from
$$
\matrix{
(i-2i')[L_{i,-2},L_{k,\ell}]\!\!\!\!&=
[[L_{i',-1},L_{i-i',-1}],L_{k,\ell}]
\vsp\hfill\cr&=
(-\a+i')[L_{i'+k,\ell-1},L_{i-i',-1}]+(-\a+i-i')[L_{i',-1},L_{i-i'+k,\ell-1}]
\vsp\hfill\cr&=
(-(-\a+i')(-\a+i-i')+(-\a+i')(-\a+i-i'))L_{i+k,-2}
\vsp\hfill\cr&=0
\mbox{ for }i'\in\Z,
\hfill\cr}
\eqno(3.70)$$
we see that $[L_{i,-2},L_{k,\ell}]=0$. If $j<-2$,
using $L_{i,j}=(-\a+i')^{-1}[L_{i',-1},L_{i-i',j+1}]$ for $i'\ne\a$,
and by induction on $-j$,
we obtain $[L_{i,j},L_{k,\ell}]=0$ for $j<-1$ and $\ell<-1$.

Assume that $j\ge1$.
If $j+\ell\ge-1$,
by writing $L_{i,j}$ as $L_{i,j}=a[L_{i',1},L_{i'',j-1}]$
for some $a\in\F,i',i''\in\Z$,
and by the first claim and induction on $j$, we have
$[L_{i,j},L_{k,\ell}]=0$. Thus assume that $j+\ell\le-2$.
Again write $L_{i,j}$ as $L_{i,j}=a[L_{i',1},L_{i'',j-1}]$
for some $a\in\F,i',i''\in\Z$. Then using induction on $j$
and the first claim, we can easily verify that
$$
[L_{i,j},L_{k,\ell}]=
[^{\ \,-\ell-2}_{-\ell-j-2}](k(j+1)-(\ell+1)i+(\ell-j)\a)L_{i+k,j+\ell},
\eqno(3.71)$$
for $k>1, l<-1$. This proves the second claim.
\par
Thus $\AA=\CC(\a)$ if ${\rm dim}(\oplus_{i\in\sZ}\AA_{i,-2})\ge2$.
\qed\par
{\bf Lemma 4.2b.} If $\b_1=-2-\b_{-1}$ and $\b_{-1}=0$ with
${\rm dim}(\oplus_{i\in\sZ}\AA_{i,-2})\le1$, then $\AA=
\BB^+(\a,-1,a_1,a_2,a'_2)$.
\par
{\bf Proof.} The result follows by a similar discussion in the last
paragraph of the proof of Lemma 4.1.
\qed\par
This completes the proof of Theorem 1.1.
\par
\vskip 10pt
\cl{\bf References}
\par
\small\baselineskip 5.7pt \lineskip 5.7pt
\parskip .08 truein
[B]  R.  Block, On torsion-free
abelian groups and Lie algebras , {\it Proc.  Amer.  Math.  Soc.},
{\bf9} (1958), 613-620.
\par
[DZ] D. Z. Djokovic and K. Zhao,
Derivations, isomorphisms, and second cohomology of generalized Block
algebras, {\it Algebra Colloquium}, {\bf3} (1996), 245-272.
\par
[KS] I.  Kaplansky and L.  J.  Santharoubane,
Harish-Chandra modules over the Virasoro algebra, MSRI Publications,
No.4 (1985), 217-231.
\par
[M]  O.  Mathieu, Classification
of simple graded Lie algebras of finite growth, {\it Invent.  Math.},
{\bf 108} (1992), 455-589.
\par
[OZ1] J. M.  Osborn and K. Zhao, Doubly
$\Z$-graded Lie algebras containing a Virasoro algebra, {\it J.  Alg.},
{\bf219} (1999), 266-298.
\par
[OZ2] J. M.  Osborn and K. Zhao, $\ZZ$-graded
Lie algebras containing a Heisenberg algebra,
{\it Nonassociative algebra and its
applications} (Sao Paulo, 1998), 259-274,
Lecture Notes in Pure and Appl. Math., 211, Dekker, New York, 2000.
\par
[OZ3] J. M.  Osborn and K. Zhao, A characterization of certain loop
algebras, {\it J.  Alg.}, {\bf221} (1999), 345-359.
\par
[OZ4] J. M.  Osborn and K. Zhao, $\ZZ$-graded Lie algebras containing
a Virasoro algebra and a Heisenberg algebra,
{\it Comm. Alg.}, {\bf29} (2001), 1677-1706.
\par
[OZ5] J. M.  Osborn and K. Zhao,
A characterization of the Block Lie algebra and its $q$-forms
in characteristic $0$,  {\it J. Alg.},
{\bf207} (1998), 367-408.
\par
[Z] K. Zhao, Isomorphic irreducible
representations of the Virasoro algebra (in Chinese), {\it Systems
Sciences and Mathematical Science}, {\bf14} (1994), 209-212.
\end{document}